# Understanding the First Order Inhomogeneous Linear Elasticity through Local Gauge Transformations

Zhihai Xiang

*Department of Engineering Mechanics, Tsinghua University, Beijing 100084, China*

Tel & Fax: +86-10-62796873

Email: xiangzhihai@tsinghua.edu.cn

**Abstract** It is well-known that classical linear elasticity equations are not form-invariant under local transformations. This is intrinsically related to the inhomogeneity of elastic media. However, the reported new linear elasticity equations for inhomogeneous media may appear in different forms. This paper tries to clarify this issue by investigating the form-invariance of the Lagrangian under local temporal or spatial gauge transformations. In this way, these new equations in different forms can be easily understood as the results from different choices of gauge fixing schemes. It recommends to choose appropriate gauges with clear physical meanings to simplify calculations.

## 1 Introduction

Although real elastic media are inhomogeneous, the classical theory of linear elasticity is always established on the assumption of homogeneity [1]. However, the effort on seeking the appropriate equations for inhomogeneous media never stops, especially when we are interested in the detailed information of the wave propagation in composites [2] or materials with micro-defects [3]. However, the corresponding theories develop very slowly.

Forty years ago, Willis established new linear elastodynamic equations for inhomogeneous media based on the variational theory and the perturbation method [4]:

$$\langle \boldsymbol{\sigma} \rangle = \boldsymbol{C}^{\text{eff}} * \langle \boldsymbol{e} \rangle + \boldsymbol{S}^{\text{eff}} \circ \langle \dot{\boldsymbol{u}} \rangle , \tag{1}$$

$$\nabla \cdot \langle \boldsymbol{\sigma} \rangle + \boldsymbol{f} = \hat{\boldsymbol{S}}^{\text{eff}} \circ \langle \dot{\boldsymbol{e}} \rangle + \boldsymbol{\rho}^{\text{eff}} \odot \langle \ddot{\boldsymbol{u}} \rangle , \tag{2}$$

where $\boldsymbol{\sigma}$ is the stress; $\boldsymbol{u}$ is the displacement; $\boldsymbol{e} = (\boldsymbol{u}\nabla + \nabla\boldsymbol{u})/2$ is the linear strain; $\boldsymbol{f}$ is the body force; the overhead dot denotes the differentiation with respect to time $t$; $\langle \ \rangle$ denotes the ensemble average operation; $\boldsymbol{C}^{\text{eff}} *$, $\boldsymbol{S}^{\text{eff}} \circ$, $\hat{\boldsymbol{S}}^{\text{eff}} \circ$, and $\boldsymbol{\rho}^{\text{eff}} \odot$ are non-local operators that involve temporal-spatial convolutions. The derivation of these equations is briefly outlined in Appendix A according to reference [2]. However, researchers had not paid so much attention to Eqs. (1) and (2), until they tried to design metamaterials [5, 6] by using transformation methods.

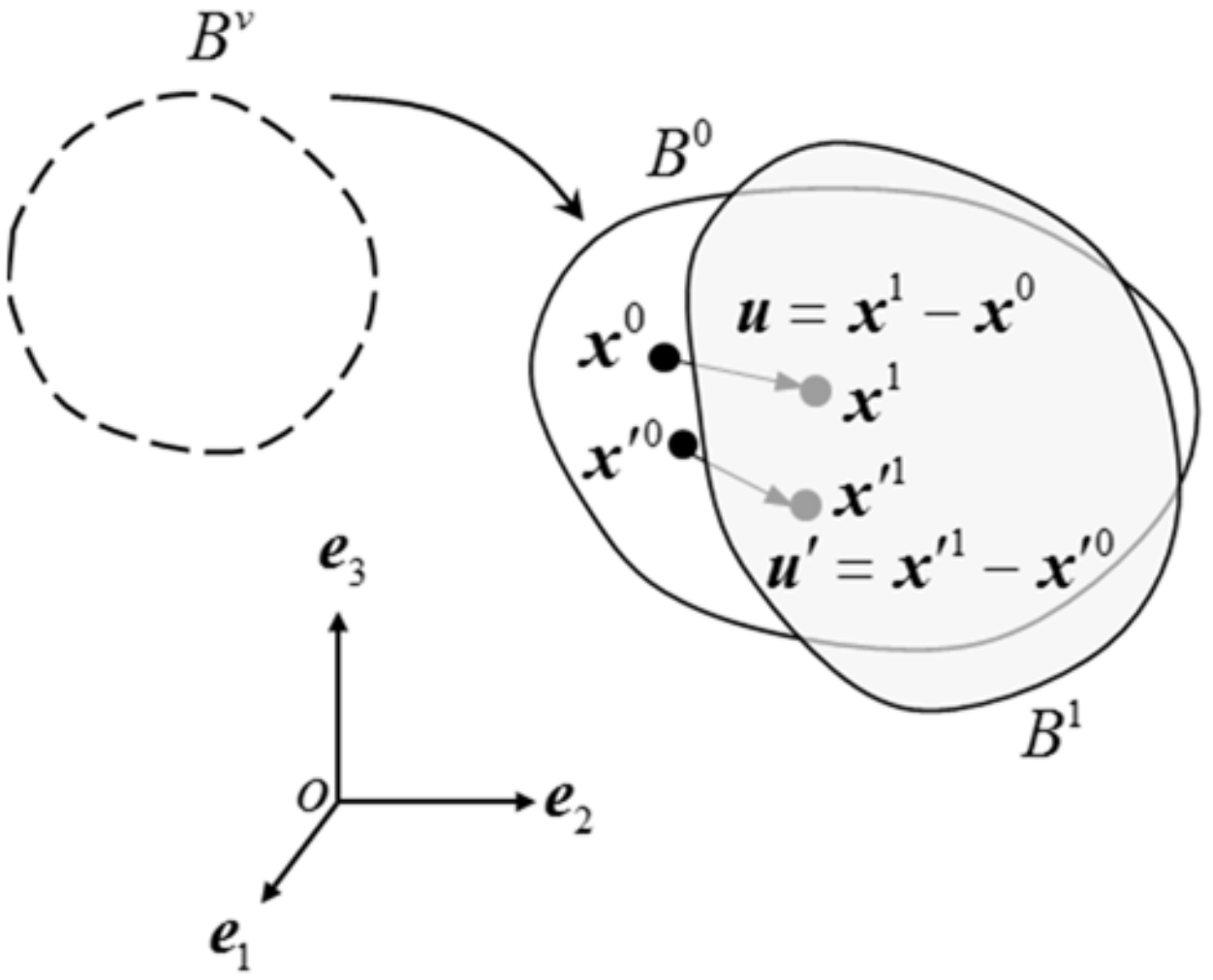

Fig. 1 The initial physical configuration $B^0$ is resulted from the transformation from the virtual configuration $B^v$. Then, $B^0$ changes to the deformed configuration $B^1$ through a small perturbation.

Metamaterials are man-made inhomogeneous media that can control wave propagation in the desired manner and achieve many exotic phenomena, such as cloaking [7], negative refraction [8], etc. An important problem in designing metamaterials is how to determine the distribution of material properties independent of the wave incident direction. To solve this problem, Pendry et al. [9] and Leonhardt [10] proposed to use the transformation method [11], by which a wave equation in the virtual homogeneous configuration $B^v$ is transformed into the physical inhomogeneous configuration $B^0$ through local mappings of coordinates and field variables (see Fig. 1). If the wave equation does not change its form after this local spatial transformation, the distribution of the material properties can be obtained by comparing the corresponding terms in the original and transformed equations. Encouraged by the great success of the transformation methods in optics [12] and acoustics [13], researchers tried to implement it to elasticity. However, Milton et al. [14] found out that the classical linear elastodynamic equations are not form-invariant after the local spatial transformation in frequency domain. At the same time, it is surprising that the transformed new equations have the similar forms as Eqs. (1) and (2). This discovery manifests the importance of Willis' original work, and then Eqs. (1) and (2) are customarily called the *Willis equations*. Since we will discuss other "Willis equations" later, Eqs. (1) and (2) are called the Classical Willis Equations (CWEs) in this paper.

Since the CWEs play non-fungible roles to design elastic metamaterials, more and more researchers are interested in investigating their properties [15-23] and even finding related applications [24]. However, there are still three questions lingering in their minds:

1) The first question is why the CWEs use velocity and strain rate coupling terms to describe lossless material (viscoelastic effects are not considered in the CWEs)? One expedient is regarding $\boldsymbol{S}^{\mathrm{eff}}$ as a complex-valued term [15-18]. Alternatively, one may suggest a coupling term related to accelerations [19-21].
2) Milton et al. [14] discuss the form-invariance under local spatial transformations in frequency domain. So, the second question is: what is the corresponding transformation in time domain that results in the CWEs? Actually, applying pure local spatial transformations to the classical elastodynamic equations in time domain, we can only obtain the following equations with displacement and strain coupling terms [25]:

$$\boldsymbol{\sigma} = \boldsymbol{C} : \boldsymbol{e} + \boldsymbol{S}^0 \cdot \boldsymbol{u}, \tag{3}$$

$$\nabla \cdot \boldsymbol{\sigma} + \boldsymbol{f} = \left(\boldsymbol{S}^0\right)^{\mathrm{T}} : \boldsymbol{e} + \boldsymbol{K}^0 \cdot \boldsymbol{u} + \boldsymbol{\rho} \cdot \ddot{\boldsymbol{u}}. \tag{4}$$

If transforming the above equations into frequency domain, we can get the same results of [14]. It should be noted that all material properties and field variables in Eqs. (3) and (4) are homogenized values of real media. Since the homogenization involves the average over a microstructure with many sub-components [26], which could move in different directions independently, we can generally obtain an anisotropic density $\boldsymbol{\rho}$ represented in tensorial form [27, 28]. This also conforms to the idea of the Jacobi metric for a general kinetic energy [29]. We have proved that $\boldsymbol{S}^0$ and $\boldsymbol{K}^0$ are related to the pre-stress $\boldsymbol{\sigma}^0$ in the initial configuration $B^0$ [30]. As Appendix B shows, this can also be understood from energy view point. Since Eqs. (3) and (4) have the similar forms as those of Eqs. (1) and (2), they are called the Willis-Form Equations (WFEs) in this paper. It has been proved that the WFEs can meet the requirement of time synchronization after applying only local spatial transformations [31]. The static version of the WFEs is not only verified by the experiment of rotational springs [32], but also implemented to identify the distribution of residual stresses [33] and even helpful to the buckling analysis of thin-walled shells [34, 35]. Although we have showed that the WFEs and the CWEs are mathematically equivalent regarding the differential property of time convolution operator [32], the underlying physical mechanism is not very clear.

3) The third question is why the CWEs involve temporal-spatial convolution operators [2, 4-6], while the "Willis equations" derived by other researchers [15-25, 30-35] use only algebraic operators?

The above questions imply that the Willis equations may present in different forms. This can be easily understood by reviewing Appendix A, in which Eqs. (A9) and (A10) just come from a single equation (A5). Willis has realized this non-uniqueness problem and tried to solve it by introducing eigen-strains [36]. He also demonstrated that the same mean stress and momentum could be obtained from equations in different forms [37]. However, it is necessary to clarify the underlying mechanism through a more general approach.

Actually, the form-invariance under transformations is not only the basis of the transformation method to design metamaterials, but also a general requirement for all physical laws. This is because coordinates are just man-made tools to facilitate calculating or measuring [38]. Therefore, the law of nature should be independent of coordinates, and we have much freedom to choose a coordinate system. It is well known that Newtonian mechanics are form-invariant under Galilean transformation, Maxwell equations are form-invariant under Lorentz transformation that coincides with the special theory of relativity. In continuum mechanics, we have the counterpart of *the principle of material frame indifference* or *objectivity* [39, 40], which requires the form-invariance under global Euclidean transformations or mappings of rigid rotations and translations. The underlying mechanism of these invariances is that a correct theory should present in the same form in any reference frame or look the same from any view point, i.e., the equation is form-invariant after coordinate transformations. This is always true not only for the aforementioned global transformations independent of coordinates, but also for local transformations that change at different temporal or spatial positions. A famous example of the form-invariance under local transformations is the general theory of relativity, which has local Lorentz invariance [41-43]. In 1918 Weyl even tried to unify electromagnetism and gravitation by investigating the invariance under the scale change of the metric tensor [44], and thus he introduced the concept of *gauge transformation* [41-46]. However, this attempt was not successful until it was adapted to the phase transformation of the wave function in quantum mechanics since 1928 [44]. The invariance of global phase transformation corresponds to the conservation of electric charge [46] and the invariance of local phase transformation coincides with electromagnetic interactions [43-45]. *The principle of gauge invariance* was further extended into the Yang-Mills theory [47] and the related Standard Model of particle physics by introducing more sophisticated transformations to describe the interactions inside atoms. This principle has become a dominant concept in theoretical physics in the second half of the 20$^{th}$ century [41-46].

Almost at the dawning of the Yang-Mills theory, Kondo [48], Bilby [49] and Kröner [50] initiated the pioneering work of connecting local gauge transformations with various deformation incompatibilities due to defects [51]. After years of developing, it has become a new branch of theory called *geometric continuum mechanics* [3, 52-56]. A related concept is *the principle of generalized objectivity*, which has been proposed to describe the form-invariance of Lagrangian under local rigid body transformations [57]. Based on this idea, a recent work studied the dynamic crack branching in brittle materials [58]. Since this kind of research is closely related to advanced differential geometry theory, strong mathematical background is required to grasp its essence. This could be a daunting challenge for many interested engineers and material scientists.

Regarding the form-invariance under transformation, we should mention the Nöether theorems [42, 46, 59-62] published in 1918, which prove the equivalence among the invariance, symmetry and conservation laws. The Nöether's first theorem talks about the conservation laws under global transformations, which has been extensively discussed in mechanics for homogeneous media [62]. The Nöether's second theorem focuses on the conservation laws under local transformations, which is closely related to inhomogeneous media and will be emphasized in this paper.

This paper tries to answer the aforementioned three questions by using local gauge theories. For this purpose, Section 2 gives a brief review of the Nöether theorems and emphasizes the form-invariance of Lagrangian under local transformations. Then, different "Willis equations" are obtained in Section 3 by using the minimal replacement method [52, 57] regarding local temporal-spatial transformations. In this way, it demonstrates that the non-uniqueness of the Willis equations comes from the gauge freedom and they can be presented in unique forms that are easy to calculate by selecting appropriate gauge fixing schemes with clear physical meanings. Finally, the answers to the three questions posed in this section are given in Section 4.

## 2 A Brief Review of the Nöether Theorems

As Fig. 1 shows, in linear elasticity we examine the deformation of a medium from the initial configuration $B^0$ to the perturbed configuration $B^1$. Correspondingly, a material point $\boldsymbol{X}$ initially located at the spatial position $\boldsymbol{x}^0$ in the global Cartesian coordinate system moves to the spatial position $\boldsymbol{x}^1$, resulting in a small displacement $\boldsymbol{u} = \boldsymbol{x}^1 - \boldsymbol{x}^0$. This deformation process should also satisfy some conservation laws, whose equivalence to the invariance or symmetry of the Hamiltonian action is revealed by the Nöether theorems [42, 46, 59-62].

Ignoring dissipation, the Hamiltonian action corresponding to the deformation from $B^0$ to $B^1$ is defined as:

$$A = \int_{\Omega} L\left(u_i, u_{i,\alpha}\right) \mathrm{d}\Omega \quad \alpha = 1,2,3,t \quad i = 1,2,3, \tag{5}$$

where $\Omega$ is a temporal-spatial domain; Roman indices only represent spatial dimensions ranging from 1 to 3; and Greek indices can also represent temporal dimension denoted by $t$. In this paper, the spatial and temporal coordinates are not coupled with each other, i.e., ignoring the relativistic effect. The comma presented in the right subscript denotes the partial differentiation with respect to a temporal-spatial coordinate $x_\alpha$.

In Eq. (5), $L$ denotes the incremental Lagrangian from $B^0$ to $B^1$:

$$L\left(u_i, u_{i,\alpha}\right) = W + \Phi - T, \tag{6}$$

where $W$, $\Phi$ and $T$ are the incremental density of strain energy, potential of external body forces and kinetic energy, respectively. In this paper, we discuss only the first order elasticity, so that $L$ is the explicit function up to the gradient of displacement and the coordinate $x_\alpha$ is implicitly included.

The stationary condition of Eq. (5) results in the Euler-Lagrange equation [61, 62]:

$$\frac{\partial L}{\partial u_i}-\left(\frac{\partial L}{\partial u_{i,\alpha}}\right)_{,\alpha}=\frac{\partial L}{\partial u_i}-\left(\frac{\partial L}{\partial u_{i,j}}\right)_{,j}-\left(\frac{\partial L}{\partial u_{i,t}}\right)_{,t}=0, \tag{7}$$

where Einstein's notation of index summation is adopted. Since $B^1$ is unknown before solving Eq. (7), $L$ should be represented in $B^0$ with known quantities in linear elasticity. In this sense, each term in Eq. (7) has its clear physical meaning:

$$\bar{f}_i=-\frac{\partial L}{\partial u_i},\ \bar{s}_{ij}=\frac{\partial L}{\partial u_{i,j}},\ \bar{p}_i=-\frac{\partial L}{\partial u_{i,t}}, \tag{8}$$

where $\bar{f}_i$ is the effective body force; $\bar{s}_{ij}$ is the total effective second Piola-Kirchhoff stress, which is usually not distinguished from the total Cauchy stress $\bar{\sigma}_{ij}$ in the classical linear elasticity; and $\bar{p}_i$ is the effective momentum. Therefore, Eq. (7) actually represents the equation of motion:

$$\bar{\sigma}_{ij,j}+\bar{f}_i=\bar{p}_{i,t}. \tag{9}$$

According to the Nöether theorems, a conservation law corresponds to a symmetry transformation of elasticity equations, or the invariance of the action $A$ after infinitesimal continuous perturbations of both coordinates and field variables:

$$x_\alpha \ \rightarrow\ x'_\alpha=x_\alpha+\mathrm{d}x_\alpha\equiv\varphi_{\alpha\beta}x_\beta, \tag{10a}$$

$$u_i(x_\alpha)\ \rightarrow\ u'_i(x'_\alpha)=u_i+\mathrm{d}u_i\equiv\psi_{ij}u_j, \tag{10b}$$

where $\varphi_{\alpha\beta}$ and $\psi_{ij}$ are the *coordinate transformation* operator and the *gauge transformation* operator, respectively. The coordinate transformation describes the change of viewpoint (or reference frame). The gauge transformation describes the corresponding change of field variables along with the coordinate transformation [43], or examines the effect of "*transferring the magnitude of a vector to an infinitesimally close point*" [44].

Since physical laws are independent of coordinates, the action $A$ should be invariant under the transformation defined in Eq. (10):

$$\int_{\Omega} L\left(u_i, u_{i,\alpha}\right) \mathrm{d}\Omega = \int_{\Omega'} L\left(u'_i, u'_{i,\alpha'}\right) \mathrm{d}\Omega' . \tag{11}$$

In this paper, "$_{,\alpha'}$" in the right subscript denotes the partial differentiation with respect to a temporal-spatial coordinate $x'_\alpha$.

According to Eq. (10), we have:

$$\mathrm{d}\Omega' = \det\left(x'_{\alpha,\beta}\right) \mathrm{d}\Omega \approx \left[1 + \left(\mathrm{d}x_\alpha\right)_{,\alpha}\right] \mathrm{d}\Omega , \tag{12a}$$

$$u'_{i,\alpha'} = \left(u_i + \mathrm{d}u_i\right)_{,\beta} x_{\beta,\alpha'} = \left(u_i + \mathrm{d}u_i\right)_{,\beta} \left[\delta_{\beta\alpha} - \left(\mathrm{d}x_\beta\right)_{,\alpha}\right] \approx u_{i,\alpha} + \left(\mathrm{d}u_i\right)_{,\alpha} - u_{i,\beta}\left(\mathrm{d}x_\beta\right)_{,\alpha} , \tag{12b}$$

where $\delta_{\alpha\beta}$ is the Kronecker delta symbol.

Substituting Eqs. (10) and (12) into Eq. (11), yields:

$$\begin{aligned} &\int_{\Omega} L\left(u_i, u_{i,\alpha}\right) \mathrm{d}\Omega \approx \int_{\Omega} L\left(u_i + \mathrm{d}u_i, u_{i,\alpha} + \left(\mathrm{d}u_i\right)_{,\alpha} - u_{i,\beta}\left(\mathrm{d}x_\beta\right)_{,\alpha}\right)\left[1 + \left(\mathrm{d}x_\alpha\right)_{,\alpha}\right] \mathrm{d}\Omega \\ &\Rightarrow \int_{\Omega} \left\{ \frac{\partial L}{\partial u_i} \mathrm{d}u_i + \frac{\partial L}{\partial u_{i,\alpha}} \left[\left(\mathrm{d}u_i\right)_{,\alpha} - u_{i,\beta}\left(\mathrm{d}x_\beta\right)_{,\alpha}\right] + L\left(\mathrm{d}x_\alpha\right)_{,\alpha} \right\} \mathrm{d}\Omega = 0 \end{aligned} . \tag{13}$$

It should be noted that Eq. (13) is "accurate" in the sense of linear elasticity by ignoring higher order small terms. The following equations also follow this rule without further repeating.

Referring to the Euler-Lagrange equation defined in Eq. (7), Eq. (13) can be rewritten as:

$$\int_{\Omega} \left[ \left( \frac{\partial L}{\partial u_{i,\alpha}} \mathrm{d}u_i \right)_{,\alpha} - \frac{\partial L}{\partial u_{i,\alpha}} u_{i,\beta} \left(\mathrm{d}x_\beta\right)_{,\alpha} + L\left(\mathrm{d}x_\alpha\right)_{,\alpha} \right] \mathrm{d}\Omega = 0 . \tag{14}$$

Since

$$\begin{aligned} -\frac{\partial L}{\partial u_{i,\alpha}} u_{i,\beta}\left(\mathrm{d}x_\beta\right)_{,\alpha} &= -\left( \frac{\partial L}{\partial u_{i,\alpha}} u_{i,\beta} \mathrm{d}x_\beta \right)_{,\alpha} + \left( \frac{\partial L}{\partial u_{i,\alpha}} \right)_{,\alpha} u_{i,\beta} \mathrm{d}x_\beta + \frac{\partial L}{\partial u_{i,\alpha}} u_{i,\alpha\beta} \mathrm{d}x_\beta \\ &= -\left( \frac{\partial L}{\partial u_{i,\alpha}} u_{i,\beta} \mathrm{d}x_\beta \right)_{,\alpha} + \left( \frac{\partial L}{\partial u_i} u_{i,\beta} + \frac{\partial L}{\partial u_{i,\alpha}} u_{i,\alpha\beta} \right) \mathrm{d}x_\beta \\ &= -\left( \frac{\partial L}{\partial u_{i,\alpha}} u_{i,\beta} \mathrm{d}x_\beta \right)_{,\alpha} + L_{,\alpha} \mathrm{d}x_\alpha \end{aligned} , \tag{15a}$$

$$L\left(\mathrm{d}x_\alpha\right)_{,\alpha} = \left(L \mathrm{d}x_\alpha\right)_{,\alpha} - L_{,\alpha} \mathrm{d}x_\alpha , \tag{15b}$$

Eq. (14) finally becomes:

$$\int_{\Omega} P_{\alpha,\alpha} \mathrm{d}\Omega = 0 , \tag{16}$$

where $P_\alpha$ is the Nöether current:

$$P_\alpha = \left( L\delta_{\alpha\beta} - \frac{\partial L}{\partial u_{i,\alpha}} u_{i,\beta} \right) \mathrm{d}x_\beta + \frac{\partial L}{\partial u_{i,\alpha}} \mathrm{d}u_i . \tag{17}$$

Since $\Omega$ is arbitrary, the integrand in Eq. (16) should be:

$$P_{\alpha,\alpha} = 0 . \tag{18}$$

This defines certain conservation laws. Especially, $L\delta_{\alpha\beta} - \partial L / \partial u_{i,\alpha} u_{i,\beta}$ corresponds to the famous Eshelby energy-momentum tensor (or the material momentum tensor) [61-64]; and $\partial L / \partial u_{i,\alpha}$ corresponds to the stress tensor (or the physical momentum tensor) [62] as indicated in Eq. (8).

If $\varphi_{\alpha\beta}$ and $\psi_{ij}$ are independent of $x_\alpha$, Eq. (10) defines global transformations, otherwise it defines local transformations. The conservation laws corresponding to global or local transformations are related to the Nöether's first or second theorem, respectively [42, 46].

For homogeneous media, the corresponding Lagrangian $L_0$ is usually easy to construct, so that researchers have tried to find out a finite number of global symmetry transformations that are related to classical conservation laws. For example, the conservation of energy is related to the global time translation; the conservation of momentum and angular momentum are related to the global physical translation and rotation, respectively; and the conservation of *J*-integral, *L*-integral and *M*-integral in fracture mechanics are related to the global material translation, rotation and scaling, respectively [62].

For inhomogeneous media, there could be an infinite number of local symmetry transformations. An important task is to construct the corresponding Lagrangian $L_1$, which is invariant to local transformations. But the related conservation laws may have obscure physical interpretations. Although the local gauge theory has achieved great success in theoretical physics [41-47], its application to mechanics is just limited in geometric continuum mechanics [3, 48-58], which regards defects in the initial configuration $B^0$ as transformed results from the virtual space $B^v$ (see Fig. 1) according to their geometric characteristics such as the Burgers vector of dislocations. We attempt to slightly generalize this concept in Section 3.

# 3 Local Gauge Transformations

## 3.1 Minimal Replacement Method

Along with the coordinate transformation, if $\psi_{ij}$ is a function of $x_\alpha$, the Lagrangian of homogeneous media will change its form after gauge transformation, i.e., $L_0\left(u'_i, u'_{i,\alpha'}\right) = L_0\left(\psi_{ij}u_j, \psi_{ij,\alpha'}u_j + \psi_{ij}u_{j,\alpha'}\right)$, so that Eq. (11) is not valid any more. In order to construct a new Lagrangian $L_1$ that is form-invariant to this local gauge transformation, we can use the minimal replacement method [52, 57] to define:

$$L_1 = L_0\left(u_i, u_{i;\alpha}\right). \tag{19}$$

The semicolon in Eq. (19) denotes the covariance derivative operator:

$$u_{i;\alpha} \equiv u_{i,\alpha} + \Gamma^j_{i\alpha}u_j, \tag{20}$$

where the connection $\Gamma^j_{i\alpha}$ describes the interaction among local freedoms in inhomogeneous media.

In Eq. (19), $L_0$ is the Lagrangian of homogeneous media, which may contain a constant pre-stress $\sigma^0_{ij}$. Since $\sigma^0_{ij} = \sigma^0_{ji}$ and the elasticity tensor satisfies $C_{ijkl} = C_{klij} = C_{ijlk} = C_{jikl}$, $L_0$ can be explicitly written as:

$$L_0\left(u_i, u_{i,\alpha}\right) = \underbrace{\sigma^0_{ij}u_{i,j} + \frac{1}{2}C_{ijkl}u_{i,j}u_{k,l}}_{W} \underbrace{-\left(\bar{f}^0_i + f_i\right)u_i}_{\Phi} \underbrace{-\frac{1}{2}\rho_{ij}u_{i,t}u_{j,t}}_{T}, \tag{21}$$

where $f_i$ is the incremental body force from $B^0$ to $B^1$; and $\bar{f}^0_i$ is the effective body force in $B^0$ that satisfies the equilibrium equation:

$$\sigma^0_{ij,j} + \bar{f}^0_i = 0. \tag{22}$$

Substituting Eq. (21) into Eq. (7) and referring to Eq. (22), we obtain a special form of Eq. (9), i.e., the well-known equation of motion for homogeneous media:

$$\sigma_{ij,j} + f_i = p_{i,t}, \tag{23}$$

where $\sigma_{ij} = \bar{\sigma}_{ij} - \sigma^0_{ij}$ is the incremental stress.

According to Eq. (19), the Lagrangian for inhomogeneous media can be constructed by replacing the differentiation operator in Eq. (21) with the covariance derivative operator defined in Eq. (20), i.e.:

$$
\begin{aligned}
L_1 &= \sigma_{ij}^0 u_{i;j} + \frac{1}{2} C_{ijkl} u_{i;j} u_{k;l} - \left(\bar{f}_i^0 + f_i\right) u_i - \frac{1}{2} \rho_{ij} u_{i;t} u_{j;t} \\
&= \sigma_{ij}^0 u_{i,j} + \sigma_{ij}^0 \Gamma_{ij}^r u_r + \frac{1}{2} C_{ijkl} \left( u_{i,j} u_{k,l} + 2 u_{i,j} \Gamma_{kl}^s u_s + \Gamma_{ij}^r u_r \Gamma_{kl}^s u_s \right), \qquad (24) \\
&\quad - \left(\bar{f}_i^0 + f_i\right) u_i - \frac{1}{2} \rho_{ij} \left( u_{i,t} u_{j,t} + 2 u_{i,t} \Gamma_{jt}^s u_s + \Gamma_{it}^r u_r \Gamma_{jt}^s u_s \right)
\end{aligned}
$$

where the conditions of $C_{ijkl} = C_{klij}$ and $\rho_{ij} = \rho_{ji}$ are used in the deduction.

It is easy to check that $L_1$ is form-invariant with arbitrary $\Gamma_{ij}^r$ and $\Gamma_{it}^r$, so that we have much freedom to construct $L_1$, which is certainly not unique. This is because we can freely choose a local coordinate system to describe the interaction among local freedoms. This also coincides with the observation that the Nöether current is always conservative if we add any divergence free vector in Eq. (17). For example, if $P_{\alpha,\alpha} = 0$, then $\left(P_\alpha + \varepsilon_{\alpha\beta\gamma} q_{\beta,\gamma}\right)_{,\alpha} = 0$ holds for any $q_\beta$, where $\varepsilon_{\alpha\beta\gamma}$ is the Levi-Civita symbol.

However, it could be dangerous if we have free rein to take any form of $L_1$, which may lack physical meanings and difficult to calculate or measure. Fortunately, as what will be discussed in the following sub-sections, $L_1$ can be constructed with slight freedoms by the understanding of the physical meanings of temporal-spatial transformations.

### 3.2 Physical Meanings of Temporal-Spatial Transformations

Mathematically, we can apply arbitrary coordinate and gauge transformations. However, it could be difficult to verify the resulted theory in experiments if the transformation lacks clear physical meanings [40]. In this section, we will discuss two types of local transformations that are easy to understand and can be verified in experiments.

Before proceeding to the detailed discussions of local temporal-spatial transformations, it should be noted that in order to describe the interaction among local freedoms in inhomogeneous media, we have to study the connection between adjacent material points in a gauge theory [43, 44]. This requirement of nonlocal theories for inhomogeneous media is different from local theories, where we trace the movement of a single material point. However in linear elasticity, all equations should be established in the reference configuration $B^0$. Therefore in the following deductions, the local temporal transformation is discussed at a fixed spatial position related with different material points due to the infinitesimal change of time; and the local spatial transformation is discussed at a fixed time related with different material points at infinitely close spatial positions.

We firstly discuss the temporal transformation at a fixed spatial position $\boldsymbol{x}^0$ in Fig. 1. Initially, two adjacent material points $\boldsymbol{X}$ and $\boldsymbol{X}'$ located at spatial positions $\boldsymbol{x}^0$ and $\boldsymbol{x}'^0$ in $B^0$, respectively. With the advance of time from $t$ to infinitely close $t'$, the material point $\boldsymbol{X}'$ moves to $\boldsymbol{x}^0$ in $B^1$ (i.e., let $\boldsymbol{x}'^1 = \boldsymbol{x}^0$ in Fig. 1), so that the displacement at $\boldsymbol{x}^0$ changes from $\boldsymbol{u}\left(\boldsymbol{x}^0, t\right) = u_i \boldsymbol{e}_i$ to $\boldsymbol{u}'\left(\boldsymbol{x}^0, t'\right)$. The difference between $\boldsymbol{u}'\left(\boldsymbol{x}^0, t'\right)$ and $\boldsymbol{u}\left(\boldsymbol{x}^0, t\right)$ contains $u_{i,t} \mathrm{d}t \boldsymbol{e}_i$ due to the change of time, and $\partial u_i / \partial X_j \mathrm{d}X_j \boldsymbol{e}_i$ due to the change of material points from $\boldsymbol{X}$ to $\boldsymbol{X}'$, where $\mathrm{d}\boldsymbol{X} = \boldsymbol{X}' - \boldsymbol{X} = \boldsymbol{x}'^0 - \boldsymbol{x}^0$. Since the change of the spatial coordinate for $\boldsymbol{X}'$ is $\mathrm{d}\boldsymbol{x} = \boldsymbol{x}^0 - \boldsymbol{x}'^0 = -\mathrm{d}\boldsymbol{X}$, $\partial u_i / \partial X_j = -\partial u_i / \partial x_j$. Then, we can write

$$\boldsymbol{u}'\left(\boldsymbol{x}^0, t'\right) = \left(u_i + u_{i,t} \mathrm{d}t - u_{i,j} \mathrm{d}X_j\right) \boldsymbol{e}_i = \left(u_i + u_{i,t} \mathrm{d}t - u_{i,j} X_{j,t} \mathrm{d}t\right) \boldsymbol{e}_i .$$

Notice that $\mathrm{d}X_j = x_j'^0 - x_j^0 = u_j^0$, where $u_j^0$ is the relative position change between these two material points from the stress-free state $B^v$ to $B^0$. Therefore, the material velocity $X_{j,t}$ actually equals to the initial relative velocity $u_{j,t}^0$ between adjacent material points. Thus, the incremental displacement at the spatial position $\boldsymbol{x}^0$ due to the infinitesimal change of time is:

$$\boldsymbol{u}'\left(\boldsymbol{x}^0, t'\right) - \boldsymbol{u}\left(\boldsymbol{x}^0, t\right) \approx \left(u_{i,t} - u_{i,j} u_{j,t}^0\right) \mathrm{d}t \boldsymbol{e}_i . \tag{25}$$

In this case, it is easy to find that $\mathrm{d}u_i\big|_{\boldsymbol{x}} = \left(u_{i,t} - u_{i,j} u_{j,t}^0\right) \mathrm{d}t$ in Eq. (10b).

Then, we discuss the spatial transformation at a fixed time $t$. As Fig. 1 shows, the displacements of the two adjacent material points $\boldsymbol{X}$ and $\boldsymbol{X}'$ can be represented as $\boldsymbol{u}\left(\boldsymbol{x}^0, t\right) = u_i \boldsymbol{g}_i$ and $\boldsymbol{u}'\left(\boldsymbol{x}'^0, t\right) = u_i' \boldsymbol{g}_i'$, respectively. $\boldsymbol{g}_i$ and $\boldsymbol{g}_i'$ are basis vectors in material manifold or Lagrange coordinate system. Without loss of generality, $\boldsymbol{g}_i$ is simply set as $\boldsymbol{e}_i$, so that $\boldsymbol{g}_i' = x_{k,i}' \boldsymbol{e}_k \approx \left(\delta_{ki} + u_{k,i}^0\right) \boldsymbol{e}_k$. The initial displacement gradient $u_{k,l}^0$ can represent the pre-strain $e_{kl}^0$ that connects to the pre-stress $\sigma_{ij}^0$ through the general Hooke's law $\sigma_{ij}^0 = C_{ijkl} u_{k,l}^0 = C_{ijkl} e_{kl}^0$. Therefore, the incremental displacement at a fixed time $t$ due to the infinitesimal change of spatial location is:

$$\boldsymbol{u}'\left(\boldsymbol{x}'^0,t\right)-\boldsymbol{u}\left(\boldsymbol{x}^0,t\right)\approx\left(u_i+u_{i,j}\mathrm{d}x_j\right)\left(\delta_{ki}+u^0_{k,ij}\mathrm{d}x_j\right)\boldsymbol{e}_k-u_i\boldsymbol{e}_i\approx\left(u_{i,j}+u^0_{i,jk}u_k\right)\mathrm{d}x_j\boldsymbol{e}_i . \quad (26)$$

In this case, it is easy to find that $\mathrm{d}u_i\big|_t=\left(u_{i,j}+u^0_{i,jk}u_k\right)\mathrm{d}x_j$ in Eq. (10b).

For temporally homogeneous media, the initial relative displacement $u^0_j$ does not change with time, so that $u^0_{j,t}=0$. For spatially homogeneous media, the initial displacement gradient $u^0_{i,j}$ is a constant over space, so that $u^0_{i,jk}=0$. Hence, Eqs. (25) and (26) represent global gauge transformations, and the transformed Lagrangian $L_0$ does not change its form.

For inhomogeneous media, $u^0_j$ is time variant and $u^0_{i,j}$ is position dependent, so that Eqs. (25) and (26) represent local gauge transformations. In these cases, $L_0$ changes its form unless we can construct a form-invariant $L_1$ according to Eq. (19). For these gauge transformations, the specific formulations of Eq. (20) can be written as:

$$u_{i;t}=u_{i,t}-u_{i,j}u^0_{j,t} , \quad (27a)$$

$$u_{i;j}=u_{i,j}+u^0_{i,jk}u_k . \quad (27b)$$

In spatial transformations, $u^0_{i,jk}$ is the connection $\Gamma^k_{ij}$ (or gauge potential) in the material manifold and $u^0_{i,j}$ serves as a gauge field. It is easy to check that $u'_{i;j'}=\left(\psi_{ik}u_k\right)_{;j'}=\psi_{ik}u_{k;j'}$ when we adopt the transformation operator $\psi_{ik}=x_{k,i'}$ as that in [14]. Notice that $u^0_{k,l}$ is closely related to the pre-stress $\sigma^0_{ij}$ through the general Hooke's law $\sigma^0_{ij}=C_{ijkl}u^0_{k,l}=C_{ijkl}e^0_{kl}$. Actually, it is reasonable to take the pre-strain or the pre-stress to measure spatial inhomogeneity, because any spatial inhomogeneity (including defects, phase transformations, plastic deformations, constrained thermal deformations, etc.) is companied with inhomogeneous pre-strains or pre-stresses.

Although the explicit formulation of the connection $\Gamma^j_{it}$ (and the corresponding $\psi_{ij}$) in temporal transformations is not clear, we can conveniently regard $-u_{i,j}u^0_{j,t}$ as $\Gamma^j_{it}u_j$, and the initial relative velocity $u^0_{j,t}$ serves as a gauge field.

In a word, $u^0_{j,t}$ and $u^0_{i,j}$ describe the temporal and spatial interactions among local freedoms, so that the covariance derivative operators defined in Eq. (27) have clear physical meanings. In addition, it is possible to measure $u^0_{j,t}$ indirectly [20, 21, 23], and pre-strains or pre-stresses can be routinely measured (e.g. by X-ray method).

The minimal replacement method has also been implemented to geometric continuum mechanics to construct the specific Lagrangian with respect to certain micro-defects, such as dislocations or disclinations [52] and damage [58], etc. In contrast to the connections constructed in this theory, Eq. (27) is just a small step forward not only using $u^0_{j,t}$ to describe temporal transformations, but also using $u^0_{i,j}$ to describe the general spatial transformations regardless of the detailed history from $B^v$ to $B^0$.

### 3.3 Local Gauge transformation for the Classical Willis Equations

Ignoring spatial transformation, a simple case of Eq. (24) is just setting $\Gamma^r_{ij} = 0$. According to Eq. (27a), and neglecting the higher order small term $\rho_{ij}\Gamma^r_{it}u_r\Gamma^s_{jt}u_s$, Eq. (24) becomes:

$$L_1 = \sigma^0_{ij}u_{i,j} + \frac{1}{2}C_{ijkl}u_{i,j}u_{k,l} - \left(\bar{f}^0_i + f_i\right)u_i - \frac{1}{2}\rho_{ij}\left(u_{i,t}u_{j,t} - 2u_{i,t}u^0_{s,t}u_{j,s}\right). \quad (28)$$

Substituting Eq. (28) into Eq. (8) and neglecting the difference between the Cauchy stress and the second Piola-Kirchhoff stress, we obtain:

$$\bar{f}_i = -\frac{\partial L_1}{\partial u_i} = \bar{f}^0_i + f_i, \quad (29)$$

$$\bar{p}_i = -\frac{\partial L_1}{\partial u_{i,t}} = \rho_{ij}u_{j,t} - \rho_{ij}u_{j,k}u^0_{k,t}, \quad (30)$$

$$\sigma_{ij} = \bar{\sigma}_{ij} - \sigma^0_{ij} = \frac{\partial L_1}{\partial u_{i,j}} - \sigma^0_{ij} = C_{ijkl}u_{k,l} + \rho_{ik}u^0_{j,t}u_{k,t}. \quad (31)$$

According to Eqs. (9) and (22) and ignoring $u^0_{k,tt}$, we get the linear equation of motion:

$$\sigma_{ij,j} + f_i = -\rho_{ij}u^0_{k,t}u_{j,kt} + \rho_{ij}u_{j,tt}. \quad (32)$$

Similar to the CWEs, the constitutive equation (31) contains a velocity coupling term and the equation of motion (32) has a term coupled with $u_{j,kt}$. However, they are different in two aspects. One difference is that the coupling terms in the CWEs involve temporal-spatial convolution operators $S^{\mathrm{eff}}_{(ij)k}$ and $\hat{S}^{\mathrm{eff}}_{i(jk)}$ defined in Eqs. (A11c) and (A11d) in Appendix A, while Eqs. (31) and (32) just have algebraic coupling terms. Another difference is that the coupling terms in Eqs. (31) and (32) do not present such symmetries as $S^{\mathrm{eff}}_{(ij)k}$ and $\hat{S}^{\mathrm{eff}}_{i(jk)}$.

As Eq. A(5) shows, the dynamic Green's function is used to calculate displacements, so that the CWEs are actually in differential-integral formulations. Since a differential formulation could have its integral counterpart such as the boundary integration equations [65], it is not difficult to understand the first difference. This understanding can also give the answer to the third question posed in Section 1.

The second difference shows the non-uniqueness of the elastic equations for inhomogeneous media. This can be easily understood from the fact that Eqs. (A9) and (A10) in Appendix A just come from a single equation (A5). This can also be regarded as the result of gauge freedoms mentioned at the end of Section 3.1. However, similar to electrodynamics [45], we can take the advantage of gauge freedoms to obtain proper equations that are easy to calculate by adopting an appropriate gauge fixing scheme. This is just like choosing an appropriate local coordinate system to facilitate calculating. For example, we can choose:

$$\rho_{ij}\Gamma^{s}_{jt}u_s = -\frac{1}{4}\left(\rho_{ij}u^0_{s,t} + \rho_{is}u^0_{j,t}\right)\left(u_{j,s} + u_{s,j}\right) = -\frac{1}{2}\left(\rho_{ij}u^0_{s,t} + \rho_{is}u^0_{j,t}\right)e_{js}. \quad (33)$$

Substituting Eq. (33) into Eq. (24), we obtain the following equations with the similar deduction:

$$\sigma_{ij} = C_{ijkl}u_{k,l} + \frac{1}{2}\left(\rho_{jk}u^0_{i,t} + \rho_{ik}u^0_{j,t}\right)u_{k,t}, \quad (34)$$

$$\sigma_{ij,j} + f_i = -\frac{1}{2}\left(\rho_{ij}u^0_{k,t} + \rho_{ik}u^0_{j,t}\right)e_{kj,t} + \rho_{ij}u_{j,tt}. \quad (35)$$

These equations can be regarded as the differential counterparts of the CWEs, because their coupling terms present the same symmetries as $S^{\mathrm{eff}}_{(ij)k}$ and $\hat{S}^{\mathrm{eff}}_{i(jk)}$. This understanding implies that the CWEs could be regarded as the results from local temporal transformations defined in Eq. (25). Thus, it answers the second question posed in Section 1. It is also noticed that these equations will not change after adding an arbitrary constant $\Delta_i$ to $u^0_i$. This means they are form-invariant under arbitrary rigid translations, just like choosing a different origin of the local coordinate system, which is the essence of gauge freedoms.

Willis has tried to solve the non-uniqueness problem by introducing eigen-strains [36], which can be regarded as a special gauge fixing scheme. He also demonstrated that the equations in different forms can result in the same response [37]. Other "Willis equations" with complex-valued coupling terms [15-18] and acceleration coupling terms [19-21] may also be understood in the sense of gauge freedoms. However, it is fair to say that the CWEs or their differential forms of Eqs. (34) and (35) are elegant due to their symmetric forms with clear physical meanings and easy to calculate.

So far, we have shown that the CWEs can be obtained by either the local temporal gauge transformation illustrated in this section or the homogenization process briefly reviewed in Appendix A. However, in order to further clarify the underlying physical mechanism, we should go deeper into the dynamic behavior of microstructures.

As aforementioned, the tensorial density comes from the homogenization of the disordered movements of sub-components inside the microstructure. In the similar logic, we notice that the total kinetic energy density of all sub-components is $T' = \sum_{M} \rho^M u^M_{j,t} u^M_{j,t} / 2$ ($M$ denotes the $M$th sub-components), while the homogenized kinetic energy density is $T = \rho_{ij} u_{i,t} u_{j,t} / 2$. Usually $T' \neq T$ due to the presence of relative velocities among sub-components. In order to compensate the lost kinetic energy density $T' - T$ in the Lagrangian, we can convert it into the effective strain energy density $\bar{W} = W - (T' - T)$. Generally, we can use the initial relative velocity $u^0_{j,t}$ to characterize the disorder pattern of a microstructure, and $u_{i,j}$ to measure the magnitude of this disorder. Therefore, $\bar{W}$ is the function of $u_{i,t}$, $u^0_{j,t}$ and $u_{i,j}$, and consequently we have the velocity coupling term in Eqs. (31) and (34). In this sense, it is clear that the velocity coupling term does not represent damping. It is also easy to understand that $u^0_{j,t}$ must be frequency dependent and vanish at zero frequency. Therefore, it is easy to obtain large coupling terms in unsymmetrical microstructures at local resonances [20-23]. However, $\rho_{ij} \Gamma^r_{it} u_r \Gamma^s_{jt} u_s$ and $u^0_{k,tt}$ could not be ignored near the resonance and the corresponding theory could be nonlinear [61]. In addition, even the microstructure is geometrically symmetric, we can still have nonzero coupling terms at certain frequencies with unsymmetrical modes. This understanding gives the answer to the first question posed in Section 1.

It is very interesting that since $u^0_{j,t}$ serves as the inhomogeneous background velocity, $\rho_{ik}u^0_{j,t}u_{k,t}$ in Eq. (31) and (34) could represent a generalized Coriolis effect. Therefore, the CWEs have intrinsic local directionality. This property has been utilized to design topological sonic crystals [66].

Referring to Eq. (25) we can substitute $\mathrm{d}u_i\big|_{\mathbf{x}} = \left(u_{i,t} - u_{i,j}u^0_{j,t}\right)\mathrm{d}t$ into Eqs. (17) and (18), and obtain:

$$\left[\left(L_1\delta_{\alpha t} - \frac{\partial L_1}{\partial u_{i,\alpha}}u_{i,t}\right)\mathrm{d}t + \frac{\partial L_1}{\partial u_{i,\alpha}}\left(u_{i,t} - u_{i,j}u^0_{j,t}\right)\mathrm{d}t\right]_{,\alpha} = 0 \\ \left(L_1\mathrm{d}t + p_i u_{i,j}u^0_{j,t}\mathrm{d}t\right)_{,t} - \left(\bar{\sigma}_{ik}u_{i,j}u^0_{j,t}\mathrm{d}t\right)_{,k} = 0 .$$

Integrating the above relation over a volume $V$, and noticing $\mathrm{d}t$ is arbitrary, we have the following conservation law:

$$\left[\int_V \left(L_1 - p_i u_{i,j}u^0_{j,t}\right)\mathrm{d}V\right]_{,t} = \int_{\partial V} \bar{\sigma}_{ik}u_{i,j}u^0_{j,t}\mathrm{d}S_k \,, \tag{36}$$

where $\partial V$ denotes the boundary of $V$, and $p_i u_{i,j}$ is the *wave momentum* [61, 62].

### 3.4 Local Gauge transformation for the Willis-Form Equations

Ignoring temporal transformation, we can have another alternative of Eq. (24) by setting $\Gamma^s_{jt} = 0$. According to Eq. (27b), we can set:

$$\Gamma^r_{ij} = u^0_{i,jr} \,. \tag{37}$$

Since $\sigma^0_{ij} = C_{ijkl}u^0_{k,l}$ , $W^0$ defined in Appendix B can be rewritten as $W^0 = W\left(x^0_i, 0\right) = \sigma^0_{ij}u^0_{i,j}/2 = C_{ijkl}u^0_{i,j}u^0_{k,l}/2$ . With this definition, $\sigma^0_{ij} = \partial W^0/\partial u^0_{i,j}$ and $C_{ijkl} = \partial^2 W^0/\partial u^0_{i,j}\partial u^0_{k,l}$. Therefore, we have the following relations:

$$\sigma^0_{ij}\Gamma^r_{ij} = \frac{\partial W^0}{\partial u^0_{i,j}}\frac{\partial u^0_{i,j}}{\partial x_r} = W^0_{,r} \,, \tag{38a}$$

$$C_{ijkl}\Gamma^s_{kl} = \frac{\partial^2 W^0}{\partial u^0_{i,j}\partial u^0_{k,l}}\frac{\partial u^0_{k,l}}{\partial x_s} = \sigma^0_{ij,s} \,, \tag{38b}$$

$$C_{ijkl}\Gamma_{ij}^{r}\Gamma_{kl}^{s} = \frac{\partial^2 W^0}{\partial u_{i,j}^0 \partial u_{k,l}^0}\frac{\partial u_{i,j}^0}{\partial x_r}\frac{\partial u_{k,l}^0}{\partial x_s} = W_{,rs}^0 . \tag{38c}$$

Substituting Eqs. (37) and (38) into Eq. (24), we immediately get the same Lagrangian as that in Appendix B, so that the corresponding constitutive equation and the equation of motion should be the same as Eq. (B6) and Eq. (B7), respectively.

In this case, the gauge field $u_{i,j}^0$ is resulted from the interaction of material points via the gauge potential $\Gamma_{ij}^k = u_{i,jk}^0$. It should be noted that $\Gamma_{ij}^k \neq \Gamma_{ji}^k$. This implies that $\Gamma_{ij}^k$ certainly contains torsion and can naturally represent the incompatible deformations induced by micro-defects such as dislocations and disclinations, etc. [3, 51-55]. In addition, $\sigma_{ij,r}^0$ is invariant after $u_{k,l}^0$ changes to $u_{k,l}^0 + \Delta_{kl}$ as long as $\left(C_{ijkl}\Delta_{kl}\right)_{,r} = 0$. This implies that there are infinite choices of $u_{k,l}^0$ without affecting the WFEs, presenting the gauge freedom under the condition of $\left(C_{ijkl}\Delta_{kl}\right)_{,r} = 0$. As illustrated in Appendix B, the inhomogeneously distributed pre-stresses result in configurational forces, which are the sources of the coupling terms in the WFEs.

In addition, referring to Eq. (26) we can substitute $\mathrm{d}u_i\big|_t = \left(u_{i,j} + u_{i,jk}^0 u_k\right)\mathrm{d}x_j$ into Eqs. (17) and (18), and obtain:

$$\begin{gathered}\left[\left(L_1\delta_{\alpha\beta} - \frac{\partial L_1}{\partial u_{i,\alpha}}u_{i,\beta}\right)\mathrm{d}x_\beta + \frac{\partial L_1}{\partial u_{i,\alpha}}\left(u_{i,j} + u_{i,jk}^0 u_k\right)\mathrm{d}x_j\right]_{,\alpha} = 0 \\ \left(L_1\mathrm{d}x_k + \bar{\sigma}_{ik}u_{i,jr}^0 u_r \mathrm{d}x_j\right)_{,k} - \left(p_i u_{i,jr}^0 u_r \mathrm{d}x_j\right)_{,t} = 0\end{gathered} .$$

Integrating the above relation over a volume $V$, and noticing $\mathrm{d}x_j$ is arbitrary, we have the following conservation law:

$$\left(\int_V p_i u_{i,jr}^0 u_r \mathrm{d}V\right)_{,t} = \int_{\partial V}\left(L_1\delta_{kj} + \bar{\sigma}_{ik}u_{i,jr}^0 u_r\right)\mathrm{d}S_k . \tag{39}$$

## 4 Concluding remarks

In linear elasticity, since the deformed configuration $B^1$ is unknown before solving the problem, all equations should be presented in the initial configuration $B^0$. In addition, it is not necessary to either distinguish $B^1$ from $B^0$, or consider the non-local effect for homogeneous media, so that we have the classical linear elasticity theory in simple forms. However for inhomogeneous media, we have to consider the subtle effects due to the change of configurations and the nonlocal interactions among material points, so that we can find those coupling terms as discussed in Section 3 and Appendixes A and B.

With the discussion in Section 3, it is clear that the first order inhomogeneous linear elasticity can be understood through the theory of local temporal-spatial transformations. Thus, the three questions posed in Section 1 can be easily answered:

1) The coupling terms related to velocity and strain rate in the CWEs have nothing to do with damping, but are the results of the material change at a fixed spatial position during the local temporal transformation or the compensation of the lost kinetic energy density during the homogenization process.
2) The CWEs are the results of local temporal transformations. The WFEs are the results of local spatial transformations.
3) The CWEs are differential-integral equations. Their pure differential counterparts are Eqs. (34) and (35).

With the help of gauge freedoms, it is also easy to understand why there are many different forms of "Willis equations" reported in literatures. Although all of them are mathematically valid, the CWEs could be the most elegant equations considering the temporal effect, and the WFEs could reveal the contribution of pre-strain gradients or pre-stress gradients. Of course, new equations will be obtained if temporal and spatial transformations are considered simultaneously.

The coupling terms in the CWEs and the WFEs are usually very weak, so that we can safely use the classical linear elasticity equations. However, these coupling terms should not be ignored at local resonances or locations with strong inhomogeneities, such as notches, crack tips or other interfaces.

Since the CWEs and the WFEs have already been partially validated by some experimental evidences [20, 21, 23, 32, 33], the classical principle of material frame indifference or objectivity could be extended to generalized versions that consider local transformations [40, 55, 57].

Weyl always used local gauges to describe the inhomogeneity due to the interaction of local freedoms in a field [46]. Although his original idea of applying gauge transformation to the scale of field variables failed in electromagnetic and gravitational fields [44], it seems that it has been revitalized in elastic field [68].

**Acknowledgements** Thanks Dr. Jing Xiao for clarifying some basic concepts of the gauge theory in theoretical physics. The related work has been supported by the National Natural Science Foundation of China [grant number 11272168, 11672144].

## Appendix A: The Classical Willis Equations

An inhomogeneous medium with elasticity tensor $\boldsymbol{C}$ and volume density $\rho$ is regarded as the perturbed medium compared to a homogeneous medium with properties of $\bar{\boldsymbol{C}}$ and $\bar{\rho}$:

$$\boldsymbol{C} = \bar{\boldsymbol{C}} + \delta\boldsymbol{C} \quad \text{and} \quad \rho = \bar{\rho} + \delta\rho \,, \tag{A1}$$

where $\delta\boldsymbol{C}$ and $\delta\rho$ are corresponding perturbations. Thus, the incremental stress $\boldsymbol{\sigma}$ and momentum $\boldsymbol{p}$ of the inhomogeneous medium are written as:

$$\boldsymbol{\sigma} = \boldsymbol{C} : \boldsymbol{e} = \left(\bar{\boldsymbol{C}} + \delta\boldsymbol{C}\right) : \boldsymbol{e} = \bar{\boldsymbol{C}} : \boldsymbol{e} + \boldsymbol{\tau} \,, \tag{A2}$$

$$\boldsymbol{p} = \rho\dot{\boldsymbol{u}} = \left(\bar{\rho} + \delta\rho\right)\dot{\boldsymbol{u}} = \bar{\rho}\dot{\boldsymbol{u}} + \boldsymbol{\pi} \,, \tag{A3}$$

where $\boldsymbol{\tau} = \delta\boldsymbol{C} : \boldsymbol{e}$ and $\boldsymbol{\pi} = \delta\rho\dot{\boldsymbol{u}}$ are polarization terms. It should be noted that we do not distinguish $\boldsymbol{\sigma}$ from the Cauchy stress and the second Piola-Kirchhoff stress in classical linear elasticity.

Substituting Eqs. (A2) and (A3) into Eq. (23), yields:

$$\nabla\cdot\left(\boldsymbol{C} : \boldsymbol{e}\right) + \boldsymbol{f} + \nabla\cdot\boldsymbol{\tau} - \dot{\boldsymbol{\pi}} = \bar{\rho}\ddot{\boldsymbol{u}} \,. \tag{A4}$$

Since $\nabla\cdot\boldsymbol{\tau} - \dot{\boldsymbol{\pi}}$ serves as an effective body force, the solution to Eq. (A4) at temporal-spatial position $\boldsymbol{x}$ can be written as:

$$\boldsymbol{u}\left(\boldsymbol{x}\right) = \bar{\boldsymbol{u}}\left(\boldsymbol{x}\right) + \int_{\Omega'} \bar{\boldsymbol{G}}\left(\boldsymbol{x}, \boldsymbol{x}'\right)\left[\nabla\cdot\boldsymbol{\tau}\left(\boldsymbol{x}'\right) - \dot{\boldsymbol{\pi}}\left(\boldsymbol{x}'\right)\right]\mathrm{d}\Omega' \,, \tag{A5}$$

where $\bar{\boldsymbol{G}}$ is the dynamic Green's function of the comparison medium, and $\bar{\boldsymbol{u}}$ is the corresponding displacement.

Applying the integration by parts on Eq. (A5) gives:

$$\boldsymbol{u} = \bar{\boldsymbol{u}} - \boldsymbol{S} \circ \boldsymbol{\tau} - \boldsymbol{M} \circ \boldsymbol{\pi} \,. \tag{A6}$$

Then, the corresponding linear strain and velocity can be obtained as:

$$\boldsymbol{e} = \bar{\boldsymbol{e}} - \boldsymbol{S}_x \circ \boldsymbol{\tau} - \boldsymbol{M}_x \circ \boldsymbol{\pi} \,, \quad \dot{\boldsymbol{u}} = \dot{\bar{\boldsymbol{u}}} - \boldsymbol{S}_t \circ \boldsymbol{\tau} - \boldsymbol{M}_t \circ \boldsymbol{\pi} \,. \tag{A7}$$

In above equations, $\boldsymbol{S}$, $\boldsymbol{S_x}$, $\boldsymbol{S_t}$, $\boldsymbol{M}$, $\boldsymbol{M_x}$, and $\boldsymbol{M_t}$ are integral operators related to $\bar{\boldsymbol{G}}$ and "∘" denotes the temporal-spatial convolution as that in Eq. (A5).

With Eq. (A6), we can obtain:

$$\boldsymbol{\tau}+\delta\boldsymbol{C}:\left(\boldsymbol{S}_x\circ\boldsymbol{\tau}+\boldsymbol{M}_x\circ\boldsymbol{\pi}\right)=\delta\boldsymbol{C}:\bar{\boldsymbol{e}},\ \ \boldsymbol{\pi}+\delta\rho\left(\boldsymbol{S}_t\circ\boldsymbol{\tau}+\boldsymbol{M}_t\circ\boldsymbol{\pi}\right)=\delta\rho\dot{\bar{\boldsymbol{u}}}. \quad \text{(A8)}$$

After applying the ensemble averaging on Eq. (A8) and eliminating $\bar{\boldsymbol{e}}$ and $\dot{\bar{\boldsymbol{u}}}$, we finally obtain the homogenized elastodynamic equations (the classical Willis Equations):

$$\langle\boldsymbol{\sigma}\rangle=\boldsymbol{C}^{\mathrm{eff}}*\langle\boldsymbol{e}\rangle+\boldsymbol{S}^{\mathrm{eff}}\circ\langle\dot{\boldsymbol{u}}\rangle, \quad \text{(A9)}$$

$$\nabla\cdot\langle\boldsymbol{\sigma}\rangle+\boldsymbol{f}=\hat{\boldsymbol{S}}^{\mathrm{eff}}\circ\langle\dot{\boldsymbol{e}}\rangle+\boldsymbol{\rho}^{\mathrm{eff}}\odot\langle\ddot{\boldsymbol{u}}\rangle, \quad \text{(A10)}$$

where

$$\boldsymbol{C}^{\mathrm{eff}}*\langle\boldsymbol{e}\rangle=\left(\bar{\boldsymbol{C}}+\langle\delta\boldsymbol{C}\rangle\right):\langle\boldsymbol{e}\rangle-\left\langle\left(\delta\boldsymbol{C}-\langle\delta\boldsymbol{C}\rangle\right):\boldsymbol{S}_x:\left(\delta\boldsymbol{C}-\langle\delta\boldsymbol{C}\rangle\right)\right\rangle\circ\langle\boldsymbol{e}\rangle, \quad \text{(A11a)}$$

$$\boldsymbol{\rho}^{\mathrm{eff}}\odot\langle\ddot{\boldsymbol{u}}\rangle=\left(\bar{\rho}+\langle\delta\rho\rangle\right)\langle\ddot{\boldsymbol{u}}\rangle-\left\langle\left(\delta\rho-\langle\delta\rho\rangle\right)\boldsymbol{M}_t\left(\delta\rho-\langle\delta\rho\rangle\right)\right\rangle\circ\langle\ddot{\boldsymbol{u}}\rangle, \quad \text{(A11b)}$$

$$\boldsymbol{S}^{\mathrm{eff}}=-\left\langle\left(\delta\boldsymbol{C}-\langle\delta\boldsymbol{C}\rangle\right):\boldsymbol{M}_x\left(\delta\rho-\langle\delta\rho\rangle\right)\right\rangle, \quad \text{(A11c)}$$

$$\hat{\boldsymbol{S}}^{\mathrm{eff}}=-\left\langle\left(\delta\rho-\langle\delta\rho\rangle\right)\boldsymbol{S}_t:\left(\delta\boldsymbol{C}-\langle\delta\boldsymbol{C}\rangle\right)\right\rangle. \quad \text{(A11d)}$$

It should be noted that $\boldsymbol{S}^{\mathrm{eff}}$ is a third-order tensor with the first two symmetric indices, i.e., $S^{\mathrm{eff}}_{(ij)k}$. $\hat{\boldsymbol{S}}^{\mathrm{eff}}$ is the adjoint of $\boldsymbol{S}^{\mathrm{eff}}$ with the last two symmetric indices, i.e., $\hat{S}^{\mathrm{eff}}_{i(jk)}$.

## Appendix B: The Willis-Form Equations

As Fig. 1 shows, if the medium in the initial configuration $B^0$ is inhomogeneous with pre-stresses $\boldsymbol{\sigma}^0$, the strain energy density in $B^0$ and the deformed configuration $B^1$ are explicit functions of the spatial position [26, 62], which are denoted as $W\left(x_i^0,0\right)$ and $W\left(x_i^1,u_{i,j}\right)$, respectively. In linear elasticity, $W\left(x_i^1,u_{i,j}\right)$ should be represented in $B^0$ by using Taylor expansion as:

$$
\begin{aligned}
W\left(x_i^1,u_{i,j}\right) &\approx W\left(x_i^0,0\right)+\left.\frac{\partial W}{\partial u_{i,j}}\right|_{x_i=x_i^0,u_{i,j}=0} u_{i,j}+\left.\frac{\partial W}{\partial x_i}\right|_{x_i=x_i^0,u_{i,j}=0} u_i \\
&+\frac{1}{2}\left.\left(\frac{\partial^2 W}{\partial u_{i,j}\partial u_{k,l}}u_{i,j}u_{k,l}+2\frac{\partial^2 W}{\partial u_{i,j}\partial x_k}u_{i,j}u_k+\frac{\partial^2 W}{\partial x_i\partial x_j}u_iu_j\right)\right|_{x_i=x_i^0,u_{i,j}=0} \\
&\equiv W\left(x_i^0,0\right)+\sigma_{ij}^0u_{i,j}+W_{,i}^0u_i+\frac{1}{2}\left(C_{ijkl}u_{i,j}u_{k,l}+2\sigma_{ij,k}^0u_{i,j}u_k+W_{,ij}^0u_iu_j\right)
\end{aligned}
\quad, \text{(B1)}
$$

where $\sigma_{ij}^0=\left.\frac{\partial W}{\partial u_{i,j}}\right|_{x_i=x_i^0,u_{i,j}=0}$, $W_{,i}^0\equiv\left.\frac{\partial W}{\partial x_i}\right|_{x_i=x_i^0,u_{i,j}=0}$ and $C_{ijkl}\equiv\left.\frac{\partial^2 W}{\partial u_{i,j}\partial u_{k,l}}\right|_{x_i=x_i^0,u_{i,j}=0}$.

To construct the Lagrangian defined in Eq. (6), we take the same *T* in Eq. (21), but *W* and *Φ* should consider inhomogeneity:

$$
\begin{aligned}
W &= W\left(x_i^1,u_{i,j}\right)-W\left(x_i^0,0\right) \\
&=\sigma_{ij}^0u_{i,j}+W_{,i}^0u_i+\frac{1}{2}\left(C_{ijkl}u_{i,j}u_{k,l}+2\sigma_{ij,k}^0u_{i,j}u_k+W_{,ij}^0u_iu_j\right)
\end{aligned}
, \quad \text{(B2)}
$$

$$
\Phi=-\left(f_i^0+\frac{1}{2}f_{i,j}^0u_j+f_i\right)u_i, \quad \text{(B3)}
$$

where $f_i^0$ is the external body force in configuration $B^0$. Since $f_i^0$ could be inhomogeneous, its gradient is included in Eq. (B3).

According to Eq. (8), the effective body force in $B^1$ but represented in $B^0$ can be obtained as:

$$
\bar{f}_i^1=-\frac{\partial L}{\partial u_i}=f_i^0+f_{i,j}^0u_j+f_i-\sigma_{kl,i}^0u_{k,l}-W_{,i}^0-W_{,ij}^0u_j. \quad \text{(B4)}
$$

These terms in Eq. (B4) have clear physical meanings. The first three terms are the contributions from the inhomogeneous external body force. The last three are due to the inhomogeneity of *W*, so that they are configurational forces according to Eshelby's definition [63, 64].

The effective body force in $B^0$ can be obtained by setting $f_i$, $u_j$ and $u_{k,l}$ to zeros in Eq. (B4):

$$
\bar{f}_i^0=f_i^0-W_{,i}^0, \quad \text{(B5)}
$$

which also accounts for the inhomogeneity of $W^0$.

Also according to Eq. (8), the incremental stress represented at position $\boldsymbol{x}^0$ can be obtained as:

$$\sigma_{ij} = \bar{\sigma}_{ij} - \sigma_{ij}^0 = \frac{\partial L}{\partial u_{i,j}} - \sigma_{ij}^0 = C_{ijkl} u_{k,l} + \sigma_{ij,k}^0 u_k . \quad \text{(B6)}$$

Based on Eqs. (22), (B5) and (B6), the equation of motion for $\sigma_{ij}$ can be obtained according to the Euler-Lagrange equation given in Eq. (7):

$$\sigma_{ij,j} + f_i = \sigma_{kl,i}^0 u_{k,l} + \sigma_{ij,jk}^0 u_k + \rho_{ij} u_{j,tt} , \quad \text{(B7)}$$

where $-\sigma_{kl,i}^0 u_{k,l} - \sigma_{ij,jk}^0 u_k$ can be regarded as a configurational force.

With Eqs. (B6) and (B7), it is clear that $S_{(ij)k}^0 = \sigma_{ij,k}^0$ and $K_{ik}^0 = \sigma_{ij,jk}^0$ in Eqs. (3) and (4). Since $\sigma_{kl}^0 = \sigma_{lk}^0$ and $C_{ijkl} = C_{ijlk}$, we can also replace $u_{k,l}$ with $e_{kl}$.

In addition, we should distinguish the incremental second Piola-Kirchhoff stress from the incremental Cauchy stress for inhomogeneous media with pre-stresses. For example, $\sigma_{ij}$ in Eq. (B6) is actually the incremental second Piola-Kirchhoff stresses, which is nonzero even under rigid translations. This is well-known in finite deformation theory [67]. However, the incremental Cauchy stress should be zero under rigid movements.